\documentclass[english,a4paper,eqno,11pt]{article}

\usepackage{amssymb}
\usepackage{latexsym}
\usepackage[english]{babel}

\font\tenbb=msbm10 at 12pt

\def\rR{\hbox{\tenbb R}}
\def\nN{\hbox{\tenbb N}}

\def\cal{\mathcal}

\font\titre cmbx10 at 18 pt

\def\gesp{\vskip1cm}
\def\esp{\vskip .6cm}
\def\pesp{\vskip .3cm}
\def\ppesp{\vskip 1mm}

\def\ni{\noindent}
\def\di{\displaystyle}

\def\Box{$\sqcap\hskip-.262truecm\sqcup$}
\newtheorem{thm}{Theorem}[section]
\newtheorem{defn}{Definition}[section]
\newtheorem{lem}{Lemma}[section]
\newtheorem{prop}{Proposition}[section]
\newtheorem{prp}{Property}[section]
\newtheorem{cor}{Corollary}[section]
\newtheorem{Exp}{Example}[section]
\newtheorem{rem}{Remark}[section]

\usepackage{graphicx,epsfig,epic} 

\psfigdriver{dvips}

\usepackage{graphicx}
\usepackage{amscd}
\usepackage{latexsym,amssymb,euscript}
\usepackage{epic,eepic}
 \usepackage{pifont}
\begin{document}

\centerline{\titre {Topological Expansion,}}
\ppesp
\centerline{\titre {Study and Applications}}

\esp
\centerline{ H\'el\`ene Porchon}
\pesp
\centerline{\small UFR 929 Math\'ematiques }
\centerline{\small Universit\' e Pierre et Marie Curie}
\centerline{\small  4 Place Jussieu, 75252 Paris Cedex 05, France}
\centerline{helene.porchon@etu.upmc.fr}

\pesp
{\small Abstract. In this paper, we introduce the notion of expanding topological space. We define the topological expansion of a topological space via local multi-homeomorphism over coproduct topology, and we prove that the coproduct family associated to any fractal family of topological spaces is expanding. In particular, we prove that the more a topological space expands, the finer the topology of its indexed states is. Using multi-homeomorphisms over associated coproduct topological spaces, we define a locally expandable topological space and we prove that a locally expandable topological space has a topological expansion. Specifically, we prove that the expanding fractal manifold is locally expandable and has a natural topological expansion.}\pesp



\section{Introduction}

The main objective of this paper is to relate the universe expansion to a topological expansion and to provide a new understanding of the space expansion via topology. However, relating topology to the continuous deformations of the space, such as its expansion, is quite impossible since topology is invariant under any continuous deformation of the space. This situation calls to circumvent the main difficulty by constructing new topological tools via fractal topology \cite{HP} that allows the detection of continuous deformations of space (such as its expansion). The fractal manifold model \cite{BF1} is an adequate lead to these new tools since it is a valid for modeling the universe expansion \cite{BF3}, and since it is naturally associated to a fractal family of topological spaces introduced in \cite{HP}.

The plan of this paper is as follow: in the preliminary part (Section 2) a general introduction to the fractal topology is presented, followed by an introduction to the fractal manifold and its own fractal topology. The Section 3 contains the main result. A definition of topological expansion via local multi-homeomorphisms over coproduct topological spaces is given in Subsection 3.1.1. We introduce the coproduct topology associated to a fractal family of topological spaces, and we prove that it is expanding in Subsection 3.1.2. In Subsection 3.1.3, we define a locally expandable topological space and we prove that a locally expandable topological space has a topological expansion. In Subsection 3.1.4, we provide examples of topological expansion. In Subsection 3.2, the fractal manifold model is presented as a fundamental example of application of these new topological tools: we prove that the expanding fractal manifold is locally expandable and then has a topological expansion. We give an overview and derive conclusions in Section 4.

\pesp
\section{Preliminary}
In this part, we give an introduction to the fractal topology that can be found in \cite{HP} as well as an introduction to the fractal manifold model that can be found in \cite{BF1}.

\subsection{Introduction to Fractal Topology}

\begin{defn}\label{Def5}
A fractal family of topological spaces is a family
\begin{equation}\label{FF}
\Big(X_n^{j_n}, {\cal F}_n^{j_n}\Big)_{{{j_n\in I_n}\atop {n\geq 0}}}
\end{equation}
where

i) for all $n\geq0$, $I_n$ is an index set such that $\hbox{Card}\ I_{n+1}>\hbox{Card}\ I_n$;

ii) for all $n\geq0$ and for all ${j_n}\in I_n$,  $\ (X_n^{j_n},{\cal F}_n^{j_n})$ is a topological space;

iii) for each $n\geq0$, the topologies ${\cal F}_n^{j_n}$ are equivalent for all ${j_n}\in I_n$;

iv) for all $n\geq0$ and for all $j_{n+1}\in I_{n+1}$, there exists a unique $j_n\in I_{n}$ such that
\begin{equation}
X_n^{j_n}\subset X_{n+1}^{j_{n+1}}\qquad  \hbox{and}\qquad {\cal F}_n^{j_n}= \Big\{O\cap X_n^{j_n}\ /\ O\in {\cal F}_{n+1}^{j_{n+1}}\Big\};
\end{equation}

v) for all $n\geq0$, for all ${j_n}\in I_n$, there exists ${j_{n+1}}\in I_{n+1}$ such that
\begin{equation}\label{Topo}
{\cal F}_n^{j_n}\subset{\cal F}_{n+1}^{j_{n+1}}\qquad\hbox{and}\qquad{\cal F}_n^{j_n}= \Big\{O\cap X_n^{j_n}\ /\ O\in {\cal F}_{n+1}^{j_{n+1}}\Big\}.
\end{equation}
\end{defn}

\begin{defn}\label{Def4}
We call fractal topological space a family of sets $\di\Big(X_n^{j_n}\Big)_{{{j_n\in I_n}\atop {n\geq 0}}}$ endowed with a  fractal topology
$\di\Big({\cal F}_n^{j_n}\Big)_{{{j_n\in I_n}\atop {n\geq 0}}}$ such that
$\di\Big(X_n^{j_n}, {\cal F}_n^{j_n}\Big)_{{{j_n\in I_n}\atop {n\geq 0}}}$ is a fractal family of topological spaces.
\end{defn}

\subsection{Introduction to Fractal Manifold}

The fractal manifold model \cite{BF1} is built using a union of topological spaces all disjoint or all the same and families of local homeomorphisms that double the local properties of the space.

\subsubsection{Fractal-Manifold}

We introduce in this part basic notions about $\delta_0$-manifold and fractal manifold that can be found in \cite{BF1} with deeper details related to the construction.

Let us consider a strictly decreasing sequence $(\varepsilon_n)_{n\geq0}$ such that $\varepsilon_n\in ]0,1[$ for all $n\geq0$.
We denote ${\cal R}_{0}=]0,\varepsilon_0[$, and  ${\cal R}_{n}=[0,\varepsilon_n[$ for all $n>0$.
For all $n\geq0$, let us consider the variables $\delta_n$ such that $\delta_n$ varies in ${\cal R}_{n}$ and $\delta_{n+1}<\delta_n$ for all $n\geq0$.

Let $f_i$, for $i=1,2,3,$ be three continuous and nowhere
differentiable functions,
defined on the interval $[a,b]\subset\rR$, with $a<b$. For $i=1,2,3$, the associated graph of $f_i$ is given by
$\Gamma_{i,0}([a,b])=\Big\lbrace(x,y)\in\rR^2/\ y=f_i(x),
\ x\in[a,b]\Big\rbrace.$
For $i=1,2,3,$ let us consider the function $f_i(x,y)={1\over 2y}\int^{x+y}_{x-y}f_i(t)dt,$ we call forward (respectively backward) mean
function of $f_i$ the function given by:
\begin{equation}\label{E0}
f_i(x+\sigma_0{\delta_0\over2},{\delta_0\over2})=\di{\sigma_0\over\delta_0}\int_x^{x+\sigma_0\delta_0}f_i(t)dt
\quad \hbox{for}\ \sigma_0=+\ (\hbox{respectively}\ \sigma_0=-),\hbox{and}\  \delta_0 \in {\cal R}_{0},
\end{equation}
and we denote by
$\Gamma_{i,\delta_0}^{\sigma_0}$ its associated graph.

We define the translation
$T_{\delta_0} :
\prod_{i=1}^{3}\Gamma_{i{\delta_0}}^{+}\times \{{\delta_0}\}
\longrightarrow \prod_{i=1}^{3}\Gamma_{i{\delta_0}}^{-}\times
\{{\delta_0}\}$ by
$$T_{\delta_0}
\Big((a_1,b_1),(a_2,b_2),(a_3,b_3)\Big)=\Big((a_1+{\delta_0},b_1),(a_2+{\delta_0},b_2),
(a_3+{\delta_0},b_3)\Big),$$
where $(a_i,b_i)\in
\Gamma_{i{\delta_0}}^{+}$, that is to say\quad  $b_i=\di
f_i(a_i+{{\delta_0} \over 2}, {{\delta_0} \over 2})={1\over
{\delta_0}}\di \int _{a_i}^{a_i+{\delta_0}}f_i(t)dt$
for   $i=1,2,3$.

\begin{defn}\label{Def0}
For $\delta_0\in{\cal R}_{0}$, let ${\cal M}_{\delta_0}$ be an
Hausdorff topological space. We say that ${\cal M}_{\delta_0}$ is an
${\delta_0}$-manifold if for every point $x \in M_{\delta_0}$,
there exist a neighborhood $\Omega _{{\delta_0}}$ of $x$ in
$M_{\delta_0}$, a map $\varphi _{\delta_0}$, and two open sets
$V^{+} _{{\delta_0}}$ of
$\prod_{i=1}^{3}\Gamma_{i{\delta_0}}^{+}\times \{{\delta_0}\} $
and $V^{-} _{{\delta_0}}$ of
$\prod_{i=1}^{3}\Gamma_{i{\delta_0}}^{-}\times \{{\delta_0}\} $
such that $\varphi _{\delta_0} : \Omega _{{\delta_0}}
\longrightarrow V^{+} _{{\delta_0}}$, and $T_{\delta_0} \circ
\varphi _{\delta_0}: \Omega _{{\delta_0}} \longrightarrow V^{-}
_{{\delta_0}} $ are two homeomorphisms.
\end{defn}

Let ${\cal M}=\bigcup_{{\delta_0}\in {\cal R }_0}{\cal M}_{{\delta_0}}$ be a union of $\delta_0$-manifolds all disjoint or all the same, where the variable
$\delta_0$ varies in ${\cal R }_0$.

\begin{defn}\label{Def1}
We call object of $\ {\cal M}=\bigcup_{{\delta_0}\in {\cal R }_0}{\cal M}_{{\delta_0}}\ $ a set $\ X=\bigcup_{{\delta_0}\in {\cal R }_0}\{x_{\delta_0}\}$, where $x_{\delta_0}\in {\cal M}_{{\delta_0}}$ for all ${{\delta_0}\in {\cal R }_0}$.
\end{defn}

Therefore an object of ${\cal M}=\bigcup_{{\delta_0}\in {\cal R }_0}{\cal M}_{{\delta_0}}$ is a family of points that has a representative element in each ${\cal M}_{{\delta_0}}$, ${{\delta_0}\in {\cal R }_0}$.

\begin{defn}\label{Def8}
We say that ${\cal M}$ admits an internal structure $\tilde{X}$ on a point $P\in {\cal M}$ if there exists a $C^0$-parametric path
\begin{equation}
\begin{array}{lll}
\tilde{X}: {\cal R }_0& \longrightarrow & \bigcup_{\delta_0\in {{\cal R }_0}} {\cal M}_{\delta_0} \\
\quad\ \  {\varepsilon} & \longmapsto & \tilde{X}(\varepsilon)\in {\cal M}_{\varepsilon} ,
\end{array}
\end{equation}
such that for all ${\varepsilon}\in {\cal R }_0$, $Range (\tilde{X}) \cap
{\cal M}_{\varepsilon}=\Big\{\tilde{X}({\varepsilon})\Big\}\ $, and there exists ${\varepsilon}'\in {\cal R }_0$ such that
$P=\tilde{X}(\varepsilon')\in {\cal M}_{\varepsilon'} $.
\end{defn}

\begin{defn}\label{Def7}
A fractal manifold is an union of Hausdorff topological spaces all disjoint or all the same ${\cal M}=\bigcup_{{\delta_0}\in {\cal R }_0}{\cal M}_{{\delta_0}}$,
which satisfies the following properties:
$\forall{\delta_0} \in {\cal R}_0$,
${\cal M}_{\delta_0}$ is a ${\delta_0}$-manifold,
 and  $\forall P\in {\cal M}$, ${\cal M}$ admits an internal structure $\tilde{X}$ on $P$ such
that there exist a neighborhood\quad
 $\Omega(Rg(\tilde{X}))=\bigcup_{{\delta_0}\in {\cal R}_0}\Omega_{\delta_0}$,
with $\Omega_{\delta_0}$ a neighborhood of $\tilde{X}({\delta_0})$ in
$M_{\delta_0}$, two open sets $V^+=\bigcup_{{\delta_0}\in {\cal
R}_0} V_{\delta_0}^+$ and $V^-=\bigcup_{{\delta_0}\in {\cal R}_0}
V_{\delta_0}^-$, where $V_{\delta_0}^\sigma$ is an open set in
$\Pi_{i=1}^3\Gamma_{i{\delta_0}}^\sigma\times\{{\delta_0}\}$ for
$\sigma=\pm$, and there exist two families of maps
$(\varphi_{\delta_0})_{{\delta_0}\in {\cal R}_0}$ and
$(T_{\delta_0}\circ\varphi_{\delta_0})_{{\delta_0}\in {\cal
R}_0}$ such that
$\varphi_{\delta_0}:\Omega_{\delta_0}\longrightarrow
V_{\delta_0}^+ $ and
$T_{\delta_0}\circ\varphi_{\delta_0}:\Omega_{\delta_0}\longrightarrow
V_{\delta_0}^-$ are homeomorphisms for all ${{\delta_0}\in {\cal
R}_0}$.
\end{defn}

\begin{defn}\label{Def6}
i) A local chart on the fractal manifold ${\cal M}$ is a triplet $(\Omega,
\varphi, T \circ \varphi)$, where $\Omega=\bigcup _{{\delta_0}
\in {\cal R}_0}\Omega_{{\delta_0}}$ is an open set of ${\cal M}$,
$\varphi= (\varphi_{\delta_0})_{{\delta_0}\in {\cal R}_0}$ is a family of homeomorphisms $\varphi_{\delta_0}$
from $\Omega_{{\delta_0}}$ to an open set $V^+_{\delta_0}$ of\quad
$\prod_{i=1}^{3}\Gamma_{i{\delta_0}}^{+}\times \{{\delta_0}\}$,
and $T\circ \varphi=(T_{\delta_0}\circ \varphi_{\delta_0})_{{\delta_0}\in {\cal R}_0}$ is a family of homeomorphisms
$T_{\delta_0}\circ \varphi_{\delta_0}$ from
$\Omega_{{\delta_0}}$ to an open set $V^-_{\delta_0}$ of\quad
$\prod_{i=1}^{3}\Gamma_{i{\delta_0}}^{-}\times \{{\delta_0}\} $
for all ${\delta_0}\in{\cal R}_0$.

ii) The coordinates of an object $Rg(\tilde{X})\subset \Omega$ related to the local
chart $(\Omega, \varphi, T \circ \varphi )$ are the coordinates of
the object $\ \varphi(Rg(\tilde{X}))\ $ in $\ \bigcup _{{\delta_0} \in {\cal
R}_0}\prod_{i=1}^{3}\Gamma_{i{\delta_0}}^{+}\times
\{{\delta_0}\}\ $, and of the object\quad $(T \circ \varphi)(Rg(\tilde{X}))$ in
$\bigcup _{{\delta_0} \in {\cal
R}_0}\prod_{i=1}^{3}\Gamma_{i{\delta_0}}^{-}\times
\{{\delta_0}\} $.
\end{defn}

\ni{\bf Notations:}\pesp

1) For all $ n\geq0$, for all $\delta_0\in {\cal R}_0,\ldots,\delta_n\in{\cal R}_n$ and $\sigma_0=\pm,\ldots,\sigma_n=\pm$, we denote by
$\di N_{\delta_0...\delta_{n}}^{\sigma_0...\sigma_{n}}$ the following set:
\begin{equation}\label{F10}
 N_{\delta_0...\delta_{n}}^{\sigma_0...\sigma_{n}}=\prod_{i=1}^3\Gamma_{i\delta_{0}...\delta_{n}}^{\sigma_0...\sigma_{n}}
\times\{\delta_{n}\}
\times...\times\{\delta_{0}\}
\end{equation}
\ni where
$\Gamma_{i\delta_{0}...\delta_{n}}^{\sigma_0...\sigma_{n}}$ represents the graph of the function:
\begin{equation}\label{mom}
F^{\sigma_0...\sigma_{n}}_{i\delta_0...\delta_{n}}(x)={\sigma_n...\sigma_0\over\delta_n...\delta_0}\int_x^{x+\sigma_n\delta_n}
\int_{t_{n-1}}^{t_{n-1}+\sigma_{n-1}\delta_{n-1}}\ldots\int_{t_{0}}^{t_0+\sigma_0\delta_0}f_i(t)dtdt_0\ldots dt_{n-1}
\end{equation}
for all $x\in [a,b]$.

2) We denote the set $\di\bigcup_{\delta_0\in{\cal R}_0}\Big(\ldots\Big(\di\bigcup_{\delta_n\in{\cal R}_n}N_{\delta_0...\delta_{n}}^{\sigma_0...\sigma_{n}}\Big)\Big)$ by
$\di\bigcup_{\delta_0\ldots\delta_n}N_{\delta_0...\delta_{n}}^{\sigma_0...\sigma_{n}}$.
\esp

Using the previous notations, the following theorem (see \cite{BF1}) illustrates the fractal manifold model:

\begin{thm}\label{Th1}
If ${\cal M}$ is a fractal manifold,
then for all $n\geq0$, and for all $k\in [2^{n},2^{n+1}-1]\cap \nN$, there exist a family of local homeomorphisms
$\varphi_k$ and a family of translations $T_k$ such that for $\sigma_j=\pm$, $j=0,1,...,n$, one has the $2^{n}$ diagrams at the $step(n)$ given by Diagram A:
\gesp

\unitlength=1.1cm
\begin{picture}(11,7)
\put(5.1,7.2){$ {\cal M}$}
\put(3.5,5.7){\tiny$\bigcup_{\delta_0} N^+_{\delta_0}$}

\put(6,5.7){\tiny$\bigcup_{\delta_0} N^-_{\delta_0}$}


\put(1.5,3.5){\tiny$\di\bigcup_{\delta_0\delta_1}
N^{++}_{\delta_0\delta_1}$}

\put(3.5,3.5){\tiny$\di\bigcup_{\delta_0\delta_1}
N^{+-}_{\delta_0\delta_1}$}

\put(5.9,3.5){\tiny$\di\bigcup_{\delta_0\delta_1}
N^{-+}_{\delta_0\delta_1}$}

\put(7.9,3.5){\tiny$\di\bigcup_{\delta_0\delta_1}
N^{--}_{\delta_0\delta_1}$}


\put(0.95,1){\tiny$\di\bigcup_{\delta_0\delta_1\delta_2}
N^{+++}_{\delta_0\delta_1\delta_2}$}

\put(2.5,-0.3){\tiny${\di\bigcup_{\delta_0\delta_1\delta_2}
N^{++-}_{\delta_0\delta_1\delta_2}}$}

\put(3,1){\tiny$\di\bigcup_{\delta_0\delta_1\delta_2}
N^{+-+}_{\delta_0\delta_1\delta_2}$}

\put(4.5,-0.3){\tiny${\di\bigcup_{\delta_0\delta_1\delta_2}
N^{+--}_{\delta_0\delta_1\delta_2}}$}

\put(5.3,1){\tiny$\di\bigcup_{\delta_0\delta_1\delta_2}
N^{-++}_{\delta_0\delta_1\delta_2}$}

\put(6.8,-0.3){\tiny${\di\bigcup_{\delta_0\delta_1\delta_2}
N^{-+-}_{\delta_0\delta_1\delta_2}}$}

\put(7.45,1){\tiny$\di\bigcup_{\delta_2\delta_1\delta_0}
N^{--+}_{\delta_0\delta_1\delta_2}$}

\put(9,-0.3){\tiny${\di\bigcup_{\delta_0\delta_1\delta_2}
N^{---}_{\delta_0\delta_1\delta_2}}$}

\put(5.2,7){\vector(-1,-1){1}}

\put(5.4,7){\vector(1,-1){1}}


\put(3.8,5.5){\vector(-1,-1){1.6}}

\put(4.1,5.5){\vector(0,-3){1.5}}

\put(6.4,5.5){\vector(0,-3){1.5}}

\put(6.8,5.5){\vector(1,-1){1.6}}


\put(2.1,3.2){\vector(-1,-3){.6}}

\put(2.1,3.2){\vector(1,-3){1}}


\put(4.1,3.2){\vector(-1,-3){.6}}

\put(4.1,3.2){\vector(1,-3){1}}


\put(6.4,3.2){\vector(-1,-3){.6}}

\put(6.4,3.2){\vector(1,-3){1}}


\put(8.5,3.2){\vector(-1,-3){.6}}

\put(8.5,3.2){\vector(1,-3){1}}
\put(4.5,5.7){\vector(1,0){1.4}}

\put(2.5,4){\vector(1,0){1.5}}

\put(6.6,4){\vector(1,0){1.5}}

\put(1.4,.5){\vector(1,0){.9}}

\put(3.5,.5){\vector(1,0){.9}}

\put(5.7,.5){\vector(1,0){.9}}

\put(8,.5){\vector(1,0){.9}}
\put(10.74,6.5){\tiny$\delta_0$}

\put(10.74,4.5){\tiny$\delta_0,\delta_1$}

\put(10.74,1.7){\tiny$\delta_0,\delta_1,\delta_2$}



 \put(4.2,6.5){\tiny$\varphi_1$}

 \put(6.2,6.5){\tiny$T_1\circ\varphi_1$}
 \put(2.5,4.8){\tiny$\varphi_2$}

 \put(4.3,4.8){\tiny$T_2\circ\varphi_2$}

 \put(6,4.8){\tiny$\varphi_3$}

 \put(7.9,4.8){\tiny$T_3\circ\varphi_3$}

\put(1.5,2.5){\tiny$\varphi_4$}

\put(2.5,2.1){\tiny$T_4\circ\varphi_4$}

\put(3.5,2.5){\tiny$\varphi_5$}

\put(4.5,2.1){\tiny$T_5\circ\varphi_5$}

\put(5.7,2.5){\tiny$\varphi_6$}

\put(6.8,2.1){\tiny$T_6\circ\varphi_6$}

\put(7.8,2.5){\tiny$\varphi_7$}

\put(8.9,2.1){\tiny$T_7\circ\varphi_7$}

\put(5.2,5.8){\tiny$T_{1}$}

\put(3.2,4.2){\tiny$T_{2}$}

\put(6.9,4.2){\tiny$T_{3}$}

\put(1.8,.2){\tiny$T_{4}$}

\put(3.8,.2){\tiny$T_{5}$}

\put(5.9,.2){\tiny$T_{6}$}

\put(8.3,.2){\tiny$T_{7}$}

\put(1,-1.2){$\vdots\vdots$}

\put(2.6,-1.2){$\vdots\vdots$}

\put(3.3,-1.2){$\vdots\vdots$}

\put(5,-1.2){$\vdots\vdots$}

\put(5.7,-1.2){$\vdots\vdots$}

\put(7,-1.2){$\vdots\vdots$}

\put(7.7,-1.2){$\vdots\vdots$}

\put(9.5,-1.2){$\vdots\vdots$}
\put(10.5,6){\line(0,1){1}}

\put(10.5,5.5){\bf $step(0)$}

\put(10.5,4){\line(0,1){1}}

\put(10.5,3.5){\bf $step(1)$}

\put(10.5,0.5){\line(0,1){2.5}}

\put(10.5,0){\bf $step(2)$}

\put(10.5,-1.2){ \vdots}

\put(2.5,-2){\footnotesize { Diagram A: Expanding diagram of the fractal manifold ${\cal M}$.}}

\thicklines
\end{picture}

\vskip3cm
\end{thm}

\subsubsection{Fractal Topology and Fractal Manifold}

\ni{\bf Notation:}

For all $n\geq0$, we denote $\Lambda_n$ the index set of cardinal $2^{n+1}$ given by:
\begin{equation}\label{DefSet}
\Lambda_n=\{\sigma_0\ldots\sigma_{n}\ /\ \sigma_0=\pm,\ldots,\sigma_n=\pm\},
\end{equation}
and we have

$\Lambda_0=\{\sigma_0\ /\ \sigma_0=\pm\}=\{+, -\}$ with cardinal $2$,

$\Lambda_1=\{\sigma_0\sigma_{1}\ /\ \sigma_0=\pm,\sigma_1=\pm\}=\{++,+-,-+,--\}$ with cardinal $2^2$,

$\Lambda_2=\{\sigma_0\sigma_{1}\sigma_{2}\ /\ \sigma_0=\pm,\sigma_1=\pm,\sigma_2=\pm\}$

$\quad=\{+++,++-,+-+,+--,-++,-+-,--+,---\}$ with cardinal $2^3$,
etc.
For all $n\geq0$, $\delta_i\in {{\cal R}_i}$ and $\sigma_i=\pm$, for $i=1,...,n$,
the set $N_{\delta_0...\delta_{n}}^{\sigma_0...\sigma_{n}}$ given by (\ref{F10}) is a  Hausdorff topological space, and
if ${\cal T}_{\delta_0...\delta_{n}}^{\sigma_0...\sigma_{n}}$ is its associated topology,
then the set $\di\bigcup _{{\delta_0}...{\delta_{n}}} N_{\delta_0...\delta_{n}}^{\sigma_0...\sigma_{n}} $ is a topological space for the  topology
called diagonal topology in \cite{BF1} and given by
\begin{equation}\label{F7}
{\cal T}_n^{\sigma_0...\sigma_{n}}= \Big\{ \Omega=\bigcup_{{\delta_0}...\delta_{n}}\Omega_{\delta_0...\delta_{n}}^{\sigma_0...\sigma_{n}}
\big/ \ \  \Omega_{\delta_0...\delta_{n}}^{\sigma_0...\sigma_{n}}\ \in  {\cal T}_{\delta_0...\delta_{n}}^{\sigma_0...\sigma_{n}}\ \
\forall {\delta_0}\in {{\cal R}_0},... ,\forall {\delta_{n}}\in {{\cal R}_n} \Big\}.
\end{equation}

The following Lemma and Theorems can be found in \cite{HP}:

\begin{lem} \label{L1}
For all $n\geq0$ and for all $\sigma_0=\pm,\ldots,\sigma_n=\pm$, the set $\di\bigcup _{{\delta_0}...{\delta_{n}}} N_{\delta_0...\delta_{n}}^{\sigma_0...\sigma_{n}}$
is identical to the subset \quad$\di\bigcup _{{\delta_0}...{\delta_{n}}} N_{\delta_0...\delta_{n}0}^{\sigma_0...\sigma_{n}\sigma_{n+1}}$
of \quad $\di\bigcup _{{\delta_0}...{\delta_{n+1}}} N_{\delta_0...\delta_{n+1}}^{\sigma_0...\sigma_{n+1}}\ $ for $\ \sigma_{n+1}=\pm$ .
\end{lem}

\begin{thm}\label{Th5}
The family of diagonal topological spaces
$\Big(\bigcup _{{\delta_0}...{\delta_{n}}} N_{\delta_0...\delta_{n}}^{j_n},\ {\cal T}_{n}^{j_n}\Big)_
{{j_n\in\Lambda_n\atop n\geq0}}$
is a fractal family of topological spaces.
\end{thm}

\begin{thm}\label{Th2}
The fractal manifold  ${\cal M}$ is locally homeomorphic via families of local homeomorphisms to the fractal topological space
$\Big(\bigcup _{{\delta_0}...{\delta_{n}}} N_{\delta_0...\delta_{n}}^{j_{n}}\Big)_{{j_n\in\Lambda_n\atop n\geq0}}$.
\end{thm}

\section{Main Results}

\subsection{Topological Expansion}

\subsubsection{Coproduct Topology and Expanding Topological Space }

We characterize an expanding topological space as follow:

\begin{defn}\label{D5}
We say that a family of topological spaces $\Big(E_n,\mathrm{T}_n\Big)_{n\geq0}$ is expanding if for all $n\geq0$ there exists a family of topological spaces $\di\Big(E_n^{j_n}, T_n^{j_n}\Big)_{{{j_n}\in I_n}}$ indexed by a set $I_n$ such that:\ppesp

i) $Card\ I_n<\ Card\ I_{n+1}$;

ii) for all $n\geq0$ and for all $j_{n+1}\in I_{n+1}$, there exists a unique $j_n\in I_{n}$ such that
\begin{equation}
 T_n^{j_n}\subset T_{n+1}^{j_{n+1}}\qquad\hbox{and}\qquad T_n^{j_n}= \Big\{O\cap X_n^{j_n}\ /\ O\in T_{n+1}^{j_{n+1}}\Big\}.
\end{equation}

iii) $(E_n,\mathrm{T}_n)_{n\geq0}$ is the topological coproduct given by
\begin{equation}
\left\{
\begin{array}{ll}
E_n =& \coprod_{j_n\in I_n} E_n^{j_n}; \\
{\mathrm{T}}_n =& \Big\{O\subset E_n\ /\ O\cap (E_n^{j_n}\times\{j_n\})\ \hbox{open set in }\ E_n^{j_n}\times\{j_n\},\quad \forall j_n\in I_n\Big\}.
\end{array}
\right.
\end{equation}
We say that the topological spaces $(E_n,\mathrm{T}_n)$, $ n>0\ $ are the indexed expansion states of the topological space  $(E_0,\mathrm{T}_0)$, and that the topological space $(E_0,\mathrm{T}_0)$ is expanding into the topological space $(E_n,\mathrm{T}_n)$ for all $n>0$.
\end{defn}

\begin{rem}
The family $\Big(E_n,\mathrm{T}_n\Big)_{n\geq0}$ represents all the transformations of the topological space $(E_0,\mathrm{T}_0)$ due to its expansion.
\end{rem}

\begin{prop}\label{P3}
Let $(E_i,T_i)_{i\in I}$ be an indexed family of topological spaces. If $E=\coprod_{i\in I}E_i$ and $\mathrm{T}$ is the coproduct topology on $E$ given by
\begin{equation}\label{F1}
 {\mathrm{T}}=\{\ O\subset E\ /\ O\cap (E_i\times\{i\})\ \hbox{open set in}\ E_i\times\{i\}\quad \forall i\in I\},
\end{equation}
then $\mathrm{T}=\Big\{ \bigcup_{i\in I}O_i\times\{i\}\ /\ \forall i\in I,\ O_i\in T_i\ \Big\}.$
\end{prop}

\ni{\bf Proof.}
We denote the set $S=\Big\{ \bigcup_{i\in I}O_i\times\{i\}\ /\ \forall i\in I,\ O_i\in T_i\ \Big\}$. We have to prove the following:

i) the set $S$ is a topology on $E$;

ii) $\mathrm{T}=S$.

i) Since $E_i$ and $\emptyset$ are in $T_i$ for all $i\in I$, then the sets $E$ and $\emptyset$ are in $S$.
Let us consider $\ \bigcup_{i\in I}O_i\times\{i\}$ and $ \bigcup_{i\in I}{U}_i\times\{i\}$ in $\mathrm{T}$. Since for all $i\in I$, $T_i$ is a topology, then $(O_i\cap{U}_i)\in T_i$, and then
$$\Big(\bigcup_{i\in I}O_i\times\{i\}\Big)\bigcap\Big(\bigcup_{i\in I}{U}_i\times\{i\}\Big)
=\bigcup_{i\in I}(O_i
\cap{U}_i)\times\{i\}$$ is in $S$.
Let us consider for all $k\in K$, $\bigcup_{i\in I}O_{i,k}\times\{i\}\in T$. We have
$$\bigcup_{k\in K}\ \bigcup_{i\in I}O_{i,k}\times\{i\}=\bigcup_{i\in I}\ \bigcup_{k\in K}O_{i,k}\times\{i\}.$$
Since for all $i\in I$, $T_i$ is a topology, then $\bigcup_{k\in K}O_{i,k}\in T_i$, and
then $\bigcup_{k\in K}\ \bigcup_{i\in I}O_{i,k}\times\{i\}\in S$, thus the set $S$ is a topology on $E$.

ii) We have to prove the following inclusions: $\mathrm{T}\subset S$ and $S\subset\mathrm{T}$.
Let us consider $O\in \mathrm{T}$. Since $O$ is a subset of $E$, then
$$O=O\cap E=O\cap \Big(\bigcup_{i\in I}E_i\times\{i\}\Big)=\bigcup_{i\in I}\Big(O\cap (E_i\times\{i\})\Big).$$
Since $O\in \mathrm{T}$, then for all $i\in I$, $O\cap (E_i\times\{i\})$ is an open set of $E_i\times \{i\}$, then there exists $O_i\in T_i$
such that $O\cap (E_i \times \{i\})=O_i\times \{i\}$. Therefore $O=\bigcup_{i\in I}O_i\times\{i\}\ \in S$.

Inversely, let us consider $O\in S$, then $O=\bigcup_{i\in I}O_i\times\{i\}$ where for all $i\in I$, $O_i\in T_i$. For all $i^*\in I$,
$$O\cap (E_{i^*}\times\{i^*\})=\Big(\bigcup_{i\in I}O_i\times\{i\}\Big)\cap \Big(E_{i^*}\times\{i^*\}\Big)$$
$$=\bigcup_{i\in I} \Big(O_i\times\{i\}\cap E_{i^*}\times\{i^*\}\Big)=\bigcup_{i\in I} \Big(O_i\cap E_{i^*}\times\{i\}\cap \{i^*\}\Big)$$
$$=O_{i^*}\cap E_{i^*}\times\{i^*\}=O_{i^*}\times\{i^*\}$$
which is an open set of $E_{i^*}\times\{i^*\}$.
Therefore $O\in \mathrm{T}$, which provides the second inclusion and ends the proof.

\rightline\Box

\begin{defn}\label{D3}
Let $M$ be a topological space and let $\ \Big(E_i,T_i\Big)_{i\in I}\ $ be an indexed family of topological spaces. We say that $M$ is locally multi-homeomorphic to the disjoint union $\ E=\coprod_{i\in I}E_i$ together with the coproduct topology $\mathrm{T}$ if for all $x\in M$ there exist $\Omega$ neighborhood of $x$ in $M$, $V=\coprod_{i\in I}V_i \in \mathrm{T}$ and a family of maps $(\psi_i)_{i \in I}$ such that for all $i \in I$, $\psi_i: \Omega\longrightarrow V_i$ is an homeomorphism.
\end{defn}

We call $\Big(\Omega,(\psi_{i})_{i\in I}\Big)$ a local multi-chart at $x\in M$. For all $y\in \Omega$, the coordinates $\Big(\psi_{i}(y)\Big)_{i\in I}$ in $E$ are the coordinates of $y$ in the multi-chart $\Big(\Omega,(\psi_{i})_{i\in I}\Big)$.

\begin{defn}\label{D4}
Let $M$ be a topological space. We say that $M$ has a topological expansion if there exists an expanding family of topological spaces $\Big(E_n, \mathrm{T}_n\Big)_{n\geq0}$ such that

i) $M$ is locally multi-homeomorphic to the topological space $(E_n, \mathrm{T}_n)$ for all $n\geq0$.

ii) $\mathrm{T}_n\subset \mathrm{T}_{n+1}$ for all $n\geq0$.
\end{defn}

\begin{defn}\label{D1}
Suppose that $(X,{\cal T})$ and $(X',{\cal T}')$ are two topological spaces.
We say that ${\cal T}'$ is finer that ${\cal T}$ if
$$\left\{
  \begin{array}{ll}
    X & \subset X' \\
    {\cal T}  &\subset {\cal T}'
  \end{array}
\right.$$
In particular if $X=X'$, we find again the classical concept of topologies comparison.
\end{defn}

\begin{prp}
Let $(X,\cal T)$ be a topological space and $A$ be a subset of $X$ endowed with the subspace topology ${\cal T}_A=\{A\cap O/ O\in{\cal T}\}$. If $A$ is open in $X$, then $\cal T$ is finer than ${\cal T}_A$.
\end{prp}

\ni{\bf Proof.}
It is obvious since ${\cal T}_A\subset \cal T$ if and only if $A$ is open in $X$.

\rightline\Box

\subsubsection{Coproduct Topology Associated to a Fractal Family}

\begin{defn}\label{D2}
Let $\di\Big(X_n^{j_n}, {\cal F}_n^{j_n}\Big)_{{{j_n\in I_n}\atop n\geq0}}$ be a fractal family of topological spaces. For all $n\geq0$ we denote by $g_n$ the map $g_n:I_{n+1}\rightarrow I_{n}$ that associates to each $j_{n+1}\in I_{n+1}$ the unique
$j_n=g_n(j_{n+1})\in I_n$ such that ${\cal F}_n^{j_n}=\{O\cap X_n^{j_n}\ /\ O\in {\cal F}_{n+1}^{j_{n+1}}\}$.
\end{defn}

\begin{prp}\label{P2}
If $\di\Big(X_n^{j_n}, {\cal F}_n^{j_n}\Big)_{{{j_n\in I_n}\atop n\geq0}}$ is a fractal family of topological spaces, then for all $n\geq0$, the map $g_n$ is a surjection on $I_{n+1}$.
\end{prp}

\ni{\bf Proof.}
The result comes from Definition \ref{Def5}, v).

\rightline\Box

\begin{lem}\label{Pr10}
If $\ \Big(X_n^{j_n}, {\cal F}_n^{j_n}\Big)_{{{j_n\in I_n}\atop {n\geq 0}}}\ $ is a fractal family of topological spaces, then
for all $n\geq0$ and for all $\ j_{n+1}\in I_{n+1}$, there exists a unique $j_n\in I_{n}$ such that $\ {\cal F}_n^{j_n}\subset{\cal F}_{n+1}^{j_{n+1}}\ $ and ${\cal F}_n^{j_n}= \Big\{O\cap X_n^{j_n}\ /\ O\in {\cal F}_{n+1}^{j_{n+1}}\Big\}$.
\end{lem}

\ni{\bf Proof.}

By Definition \ref{Def5},iv), for all ${n\geq0}$ and for all $j_{n+1}\in I_{n+1}$, there exists a unique $j_n^*\in I_n$ such that
$$\left\{
    \begin{array}{ll}
      {\cal F}_n^{j_n^*}=\{O\cap X_n^{j_n^*}\ /\ O\in {\cal F}_{n+1}^{j_{n+1}}\} \\
      X_n^{j_n^*}\subset X_{n+1}^{j_{n+1}}
    \end{array}
  \right.
$$
By Definition \ref{Def5},v), for $j_n^*\in I_n$ there exists $j_{n+1}^*\in I_{n+1}$ such that
$$\left\{
    \begin{array}{ll}
      {\cal F}_n^{j_n^*}=\{O\cap X_n^{j_n^*}\ /\ O\in {\cal F}_{n+1}^{j_{n+1}^*}\} \\
      {\cal F}_n^{j_n^*}\subset {\cal F}_{n+1}^{j_{n+1}^*}
    \end{array}
  \right.
$$
which implies $X_n^{j_n^*}\subset X_{n+1}^{j_{n+1}^*}$. Since $X_n^{j_n^*}\subset X_{n+1}^{j_{n+1}}$ and $X_n^{j_n^*}\subset X_{n+1}^{j_{n+1}^*}$, then there exist two canonical injections $i:X_n^{j_n^*}\longrightarrow X_{n+1}^{j_{n+1}}$ and $i^*:X_n^{j_n^*}\longrightarrow X_{n+1}^{j_{n+1}^*}$ such that $i(X_n^{j_n^*})=X_n^{j_n^*}$ and $i^*(X_n^{j_n^*})=X_n^{j_n^*}$, which means that the space $X_n^{j_n^*}$ is seen in both spaces $X_{n+1}^{j_{n+1}}$ and $X_{n+1}^{j_{n+1}^*}$. Since ${\cal F}_n^{j_n^*}$ is induced by ${\cal F}_{n+1}^{j_{n+1}^*}$ and ${\cal F}_n^{j_n^*}$ is induced by ${\cal F}_{n+1}^{j_{n+1}}$, then the injections $i$ and $i^*$ are continuous. By Definition \ref{Def5}, iii), the topologies ${\cal F}_{n+1}^{j_{n+1}^*}$ and ${\cal F}_{n+1}^{j_{n+1}}$ are equivalent, then there exists a homeomorphism $\varphi: X_{n+1}^{j_{n+1}}\longrightarrow X_{n+1}^{j_{n+1}^*}$ such that
$i^*=\varphi\circ i$. We have $i^*(X_n^{j_n^*})=\varphi\circ i(X_n^{j_n^*})$, which implies that $X_n^{j_n^*}=\varphi(X_n^{j_n^*})$. Therefore if $X_n^{j_n^*}$ is an open set in $X_{n+1}^{j_{n+1}^*}$, then $\varphi^{-1}(X_n^{j_n^*})=X_n^{j_n^*}$ is an open set in $X_{n+1}^{j_{n+1}}$.
Since ${\cal F}_n^{j_n^*}$ is induced by ${\cal F}_{n+1}^{j_{n+1}}$ and $X_n^{j_n^*}$ is an open set in $X_{n+1}^{j_{n+1}}$, then ${\cal F}_n^{j_n^*}\subset {\cal F}_{n+1}^{j_{n+1}}$, which completes the proof.

\rightline\Box

\begin{defn}
Let $\di\Big(X_n^{j_n}, {\cal F}_n^{j_n}\Big)_{{{j_n\in I_n}\atop n\geq0}}$ be a fractal family of topological spaces. We define for all $n\geq0$ the coproduct $(X_n,F_n)$ associated to  $\di\Big(X_n^{j_n}, {\cal F}_n^{j_n}\Big)_{{{j_n\in I_n}\atop n\geq0}}$ to be the disjoint union
\begin{equation}\label{FO5}
X_n=\coprod_{j_n\in I_n}X_n^{j_n}=\bigcup_{j_n\in I_n}X_n^{j_n}\times\{j_n\}
\end{equation}
together with the coproduct topology ${\cal F}_n$ given by
\begin{equation}\label{FO3}
{\cal F}_n=\Big\{O\subset X_n \ /\  O\cap (X_n^{j_n}\times\{j_n\})\ \hbox{open set in}\ X_n^{j_n}\times\{j_n\}\quad \forall j_n\in I_n\Big\}.
\end{equation}
The topological spaces $\di\Big(X_n^{j_n}, {\cal F}_n^{j_n}\Big)$, $j_n\in I_n$ are called the fractal constituent spaces of the coproduct $(X_n,{\cal F}_n)$, and $\Big(X_n,{\cal F}_n\Big)_{n\geq0}$ is called the coproduct family associated to  $\di\Big(X_n^{j_n}, {\cal F}_n^{j_n}\Big)_{{{j_n\in I_n}\atop n\geq0}}$.
\end{defn}

\begin{rem}
By Proposition \ref{P3},
\begin{equation}\label{FO11}
{\cal F}_n  =\Big\{\ \bigcup_{j_n\in I_n}O_n^{j_n}\times\{j_n\}\ /\ \forall j_n\in I_n,\ O_n^{j_n}\in{\cal F}_n^{j_n}\ \Big\}.
\end{equation}
\end{rem}

\begin{thm}\label{P5}
If $\di\Big(X_n^{j_n}, {\cal F}_n^{j_n}\Big)_{{{j_n}\atop n}{\in\atop\geq}{{I_n}\atop 0}}$ is a fractal family of topological spaces, then the associated coproduct family $\Big(X_n, {\cal F}_n\Big)_{n\geq0}$ is expanding.
\end{thm}

\ni{\bf Proof.}
Let  $\di\Big(X_n^{j_n}, {\cal F}_n^{j_n}\Big)_{{{j_n\in I_n}\atop n\geq0}}$ be a fractal family of topological spaces, and let $\Big(X_n, {\cal F}_n\Big)_{n\geq0}$ be its associated coproduct family.

We have to verify that the conditions of the Definition \ref{D5} are satisfied. By Definition \ref{Def5}, $Card\ I_n<Card\ I_{n+1}$.
By Lemma \ref{Pr10}, for all ${n\geq0}$ and for all $j_{n+1}\in I_{n+1}$, there exists a unique $j_n\in I_n$ such that
$$\left\{
    \begin{array}{ll}
      {\cal F}_n^{j_n}=\{O\cap X_n^{j_n}\ /\ O\in {\cal F}_{n+1}^{j_{n+1}}\} \\
      {\cal F}_n^{j_n}\subset {\cal F}_{n+1}^{j_{n+1}}
    \end{array}
  \right.
$$
and for all $n\geq0$, the space $X_n=\coprod_{j_n\in I_n}X_n^{j_n}$ is endowed with the coproduct topology ${\cal F}_n$, which completes the proof.

\rightline\Box

\begin{cor}
If $\di\Big(X_n^{j_n}, {\cal F}_n^{j_n}\Big)_{{{j_n}\atop n}{\in\atop\geq}{{I_n}\atop 0}}$ is a fractal family of topological spaces and  $\Big(X_n, {\cal F}_n\Big)_{n\geq0}$ its associated coproduct family, then the topological space $(X_0, {\cal F}_0)$ is expanding into the topological space $(X_n, {\cal F}_n)$ for all $n\geq0$.
\end{cor}

\ni{\bf Proof.}
The proof is straight forward from Theorem \ref{P5} and Definition \ref{D5}.

\rightline\Box

\begin{lem}\label{Pr8}
Let $\di\Big(X_n^{j_n}, {\cal F}_n^{j_n}\Big)_{{{j_n\in I_n}\atop {n\geq 0}}}$ be a fractal family of topological spaces.
For all $n\geq0$, for all $j_n\in I_n$ and for all $j_{n+1}\in I_{n+1}$ such that $g_n(j_{n+1})=j_n$, the space $X_n^{j_n}\times\{j_n\}$ is identified to the subspace $X_{n}^{j_n}\times\{j_{n+1}\}$ of $X_{n+1}^{j_{n+1}}\times\{j_{n+1}\}$.
\end{lem}

\ni{\bf Proof.}
Let $n\geq0$ and $j_n\in I_n$. Let $j_{n+1}\in I_{n+1}$ such that $g_n(j_{n+1})=j_n$. Since $g_n(j_{n+1})=j_n$, by Lemma \ref{Pr10}, we have ${\cal F}_n^{j_n}\subset {\cal F}_{n+1}^{j_{n+1}}$ and ${\cal F}_n^{j_n}= \Big\{O\cap X_n^{j_n}\ /\ O\in {\cal F}_{n+1}^{j_{n+1}}\Big\}$. Let us consider the following canonical injection
\begin{equation}
\begin{array}{lll}
\zeta: X_n^{j_n}\times\{j_n\}& \longrightarrow & X_{n+1}^{j_{n+1}}\times\{j_{n+1}\}\\
\qquad\qquad  (x,j_n) & \longmapsto & (x,j_{n+1})
\end{array}
\end{equation}
The injection $\zeta$ is continuous: indeed let $O=O_{n+1}^{j_{n+1}}\times\{j_{n+1}\}$ be an open set of $X_{n+1}^{j_{n+1}}\times\{j_{n+1}\}$, where $O_{n+1}^{j_{n+1}}\in {\cal F}_{n+1}^{j_{n+1}}$, then
$$\zeta^{-1}(O)=\zeta^{-1}(O_{n+1}^{j_{n+1}}\times\{j_{n+1}\})=(O_{n+1}^{j_{n+1}}\cap X_n^{j_n})\times\{j_n\}.$$
Since ${\cal F}_n^{j_n}$ is induced by ${\cal F}_{n+1}^{j_{n+1}}$, then $O_{n+1}^{j_{n+1}}\cap X_n^{j_n}\in {\cal F}_n^{j_n}$, then $\zeta^{-1}(O)$ is an open set of $X_n^{j_n}\times\{j_n\}$.
Therefore the continuous injection $\zeta$ is an homeomorphism from $X_n^{j_n}\times\{j_n\}\ $ to $\ \zeta\Big(X_n^{j_n}\times\{j_n\}\Big)=X_n^{j_n}\times\{j_{n+1}\}$, which allows to identify the space $X_n^{j_n}\times\{j_n\}\ $ to the space $\ X_n^{j_n}\times\{j_{n+1}\}$.

\rightline\Box

\begin{lem}\label{Pr9}
If $\di\Big(X_n^{j_n}, {\cal F}_n^{j_n}\Big)_{{{j_n\in I_n}\atop {n\geq 0}}}$ is a fractal family of topological spaces, then for $n\geq0$
\begin{equation}
\di\bigcup_{j_{n+1}\in I_{n+1}}O_{n+1}^{j_{n+1}}\times\{j_{n+1}\}=\bigcup_{j_{n}\in I_{n}}\Big(\bigcup_{{j_{n+1}\in I_{n+1}\atop g_n(j_{n+1})=j_n}}O_{n+1}^{j_{n+1}}\times\{j_{n+1}\}\Big),
\end{equation}
where $O_{n+1}^{j_{n+1}}\in {\cal F}_{n+1}^{j_{n+1}}$ for all $j_{n+1}\in I_{n+1}$.
\end{lem}

\ni{\bf Proof.} Let $n\geq0$.
Since $Card\ I_n<Card\ I_{n+1}$, the map $g_n$ is surjective but not injective. For a given $j_n\in I_n$, we consider the set of all $j_{n+1}\in I_{n+1}$ that have the same image $j_n$ by the map $g_n$, that is to say such that $g_n(j_{n+1})=j_n$. Therefore we can write the union $\di\bigcup_{j_{n+1}\in I_{n+1}}O_{n+1}^{j_{n+1}}\times\{j_{n+1}\}\ $ as $\di\ \bigcup_{j_{n}\in I_{n}}\Big(\bigcup_{{j_{n+1}\in I_{n+1}\atop g_n(j_{n+1})=j_n}}O_{n+1}^{j_{n+1}}\times\{j_{n+1}\}\Big)$, which gives the property.

\rightline\Box

\begin{rem}
Another proof of the Lemma \ref{Pr9} can be obtained by proving the following double inclusion:
$$
\left\{
  \begin{array}{ll}
    \di\bigcup_{j_{n+1}\in I_{n+1}}O_{n+1}^{j_{n+1}}\times\{j_{n+1}\}\subset \bigcup_{j_{n}\in I_{n}}\Big(\bigcup_{{j_{n+1}\in I_{n+1}\atop g_n(j_{n+1})=j_n}}O_{n+1}^{j_{n+1}}\times\{j_{n+1}\}\Big), &  \\
    \di\bigcup_{j_{n}\in I_{n}}\Big(\bigcup_{{j_{n+1}\in I_{n+1}\atop g_n(j_{n+1})=j_n}}O_{n+1}^{j_{n+1}}\times\{j_{n+1}\}\Big)\subset \di\bigcup_{j_{n+1}\in I_{n+1}}O_{n+1}^{j_{n+1}}\times\{j_{n+1}\}. &
  \end{array}
\right.
$$
\end{rem}

\begin{thm}\label{P1}
If $\di\Big(X_n^{j_n}, {\cal F}_n^{j_n}\Big)_{{{j_n\in I_n}\atop {n\geq 0}}}$ is a fractal family of topological spaces, and $\Big(X_n, {\cal F}_n\Big)_{n\geq0}$ is its associated coproduct family, then
\begin{equation}
\forall {n\geq0},\quad {\cal F}_n\subset {\cal F}_{n+1}.
\end{equation}
\end{thm}

\ni{\bf Proof.}
Let  $\di\Big(X_n^{j_n}, {\cal F}_n^{j_n}\Big)_{{{j_n\in I_n}\atop {n\geq 0}}}$ be a fractal family of topological spaces and let ${n\geq0}$.
Let us consider $\Omega\in  {\cal F}_n$. By (\ref{FO11}), $\Omega=\di\bigcup_{j_n\in I_n}\Omega_n^{j_n}\times\{j_n\}$, where $\Omega_n^{j_n}\in {\cal F}_n^{j_n}$ for all $j_n\in I_n$.
By Definition \ref{Def5}, v), for all $j_n\in I_n$, there exists $j_{n+1}\in I_{n+1}$ such that ${\cal F}_n^{j_n}\subset {\cal F}_{n+1}^{j_{n+1}}$.
Let us denote $p=Card\ \{j_{n+1}\in I_{n}\ /\ g(j_{n+1})=j_n \}$.
Since $g_n$ is surjective but not injective, then $p\geq1$.
We denote by $j_{n+1}^1,..., j_{n+1}^{p}$ all the elements of the set
$\{j_{n+1}\in I_{n}\ /\ g(j_{n+1})=j_n \}$, and we have
$$\left\{
  \begin{array}{ll}
    {\cal F}_n^{j_n}\subset & {\cal F}_{n+1}^{j_{n+1}^1}; \\
    \ldots & \ldots \\
    {\cal F}_n^{j_n}\subset & {\cal F}_{n+1}^{j_{n+1}^p},
  \end{array}
\right.
$$
then there exist $O_{n+1}^{j_{n+1}^1}\in {\cal F}_{n+1}^{j_{n+1}^1}$,$\ldots$ , $O_{n+1}^{j_{n+1}^p}\in {\cal F}_{n+1}^{j_{n+1}^p}$ such that
\begin{equation}
\Omega_n^{j_n}=O_{n+1}^{j_{n+1}^1}=\ldots=O_{n+1}^{j_{n+1}^p},
\end{equation}
then we can write
$$\di\Omega_n^{j_n}=\bigcup_{{j_{n+1}\in I_{n+1}}\atop {g_n(j_{n+1})=j_n}}O_{n+1}^{j_{n+1}}$$
and we have
$$\Omega_n^{j_n}\times\{j_n\}=\Big(\bigcup_{{j_{n+1}\in I_{n+1}}\atop {g_n(j_{n+1})=j_n}}O_{n+1}^{j_{n+1}}\Big)\times\{j_n\}$$
$$=\Big(O_{n+1}^{j_{n+1}^1}\cup\ldots\cup O_{n+1}^{j_{n+1}^p}\Big)\times\{j_n\}$$
$$=\bigcup_{{j_{n+1}\in I_{n+1}}\atop {g_n(j_{n+1})=j_n}}\Big(O_{n+1}^{j_{n+1}}\times\{j_n\}\Big).$$
Since for all $j_{n+1}\in I_{n+1}$ such that $g_n(j_{n+1})=j_n$, we have
$O_{n+1}^{j_{n+1}}\subset X_n^{j_n}$, then by Lemma \ref{Pr8}, the set $O_{n+1}^{j_{n+1}}\times\{j_n\}$ can be identified to the set $O_{n+1}^{j_{n+1}}\times\{j_{n+1}\}$, and we obtain
$$\Omega_n^{j_n}\times\{j_n\}=\bigcup_{{j_{n+1}\in I_{n+1}}\atop {g_n(j_{n+1})=j_n}}\Big(O_{n+1}^{j_{n+1}}\times\{j_{n+1}\}\Big),$$
therefore
$$\Omega=\bigcup_{j_n\in I_n}\Omega_n^{j_n}\times\{j_n\}=\bigcup_{j_n\in I_n}\Big(\bigcup_{{j_{n+1}\in I_{n+1}}\atop {g_n(j_{n+1})=j_n}}\Big(O_{n+1}^{j_{n+1}}\times\{j_{n+1}\}\Big)\Big).$$
By Lemma \ref{Pr9},
$\ \Omega=\di\bigcup_{j_{n+1}\in I_{n+1}}O_{n+1}^{j_{n+1}}\times\{j_{n+1}\}\ \in {\cal F}_{n+1},$
which ends the proof.

\rightline\Box

\begin{cor}\label{P4}
If $\di\Big(X_n^{j_n}, {\cal F}_n^{j_n}\Big)_{{{j_n}\atop n}{\in\atop\geq}{{I_n}\atop 0}}$ is a fractal family of topological spaces, and $\Big(X_n, {\cal F}_n\Big)_{n\geq0}$ is its associated coproduct family, then  $\di\Big(X_n\Big)_{n\geq 0}$ is ascending.
\end{cor}

\ni{\bf Proof.}
Its a direct consequence of Theorem \ref{P1}.

\rightline\Box

\begin{cor}
If $\di\Big(X_n^{j_n}, {\cal F}_n^{j_n}\Big)_{{{j_n\in I_n}\atop {n\geq 0}}}$ is a fractal family of topological spaces, and $(X_n,{\cal F}_n)_{n\geq0}$ is its associated coproduct family, then the more the topological space $(X_0, {\cal F}_0)$ expands, the finer the topology of its indexed states is.
\end{cor}

\ni{\bf Proof.}
For all $n\geq0$, we have $X_n\subset X_{n+1}$ by Corollary \ref{P4} and ${\cal F}_n\subset{\cal F}_{n+1}$ by Theorem \ref{P1}, then by Definition \ref{D1} we conclude that ${\cal F}_{n+1}$ is finer than ${\cal F}_n$, which gives the result.

\rightline\Box

\subsubsection{Locally Expandable Topological Space}

\begin{defn}\label{D6}
A topological space $M$ is said to be locally expandable if there exists a fractal family of topological spaces $\di\Big(X_n^{j_n}, {\cal F}_n^{j_n}\Big)_{{{j_n\in I_n}\atop {n\geq 0}}}$ such that $M$ is locally multi-homeomorphic to the associated coproduct topological space $(X_n,{\cal F}_n)$ for all $n\geq0$.
\end{defn}

\begin{Exp}
The fractal manifold introduced in \cite{BF1} is an example of locally expandable topological space.
\end{Exp}

\begin{thm}\label{Th3}
If $M$ is a locally expandable topological space, then $M$ has a topological expansion.
\end{thm}

\ni{\bf Proof.}
Let $M$ be a locally expandable topological space. By Definition \ref{D6}, there exists a fractal family of topological spaces $\di\Big(X_n^{j_n}, {\cal F}_n^{j_n}\Big)_{{{j_n\in I_n}\atop {n\geq 0}}}$ such that $M$ is locally multi-homeomorphic to $(X_n, {\cal F}_n)$ for all $n\geq0$. By Theorem \ref{P5}, the family of topological spaces $\Big(X_n, {\cal F}_n\Big)_{n\geq0}$ is expanding, then there exists an expanding family of topological spaces $\Big(X_n, {\cal F}_n\Big)_{n\geq0}$ such that $M$ is locally multi-homeomorphic to the associated coproduct $(X_n, {\cal F}_n)$ for all $n\geq0$. Moreover using Theorem \ref{P1}, we have ${\cal F}_n\subset {\cal F}_{n+1}$ for all $n\geq0$, therefore by Definition \ref{D4}, $M$ has a topological expansion.

\rightline\Box

\subsubsection{Examples of Topological Expansion by Stretching}\label{Ex}

To give an example of topological space that has a topological expansion, we have first to build a fractal family of topological spaces.\pesp
\begin{description}
  \item[Stretching Process on $\rR$:] Let $]a,b[$ be an open interval of $\rR$ with $a<b$. We define a sequence of open intervals $(X_n^{j_n})_{j_n\in I_n}$ of $\rR$ endowed with the induced topology of $\rR$  by stretching the interval $]a,b[$ left and right step by step as follow: let us consider a strictly decreasing sequence of positive real numbers $(\varepsilon_n)_{n\geq0}$.
\end{description}
 \begin{itemize}
      \item For the $step(0)$, by subtracting $\varepsilon_0$ to the left of the interval $]a,b[$, we obtain $X_0^-=]a-\varepsilon_0,b[$, and by adding $\varepsilon_0$ to the right of the interval $]a,b[$, we obtain $X_0^+=]a,b+\varepsilon_0[$;
      \item For the $step(1)$, by subtracting $\varepsilon_1$ to the left of the interval $X_0^-$ and by adding $\varepsilon_1$ to the right of the interval $X_0^-$, we obtain the open intervals $$X_1^{--}=]a-\varepsilon_0-\varepsilon_1,b[\quad \hbox{and}\quad X_1^{-+}=]a-\varepsilon_0,b+\varepsilon_1[,$$ by subtracting $\varepsilon_1$ to the left of the interval $X_0^+$ and by adding $\varepsilon_1$ to the right of the interval $X_0^+$, we obtain the intervals $$X_1^{+-}=]a-\varepsilon_1,b+\varepsilon_0[\quad \hbox{and}\quad X_1^{++}=]a,b+\varepsilon_0+\varepsilon_1[;$$
\item For the $step(2)$, we stretch each interval of the previous step left and right by subtracting $\varepsilon_2$ for the left stretching and adding $\varepsilon_2$ for the right stretching, to obtain:
 $$X_2^{---}=]a-\varepsilon_0-\varepsilon_1-\varepsilon_2,b[, \quad \hbox{and}\quad X_2^{--+}=]a-\varepsilon_0-\varepsilon_1,b+\varepsilon_2[$$
$$X_2^{-+-}=]a-\varepsilon_0-\varepsilon_2,b+\varepsilon_1[\quad \hbox{and}\quad X_2^{-++}=]a-\varepsilon_0,b+\varepsilon_1+\varepsilon_2[$$
$$X_2^{+--}=]a-\varepsilon_1-\varepsilon_2,b+\varepsilon_0[\quad \hbox{and}\quad X_2^{+-+}=]a-\varepsilon_1,b+\varepsilon_0+\varepsilon_2[$$
$$X_2^{++-}=]a-\varepsilon_2,b+\varepsilon_0+\varepsilon_1[\quad \hbox{and}\quad X_2^{+++}=]a,b+\varepsilon_0+\varepsilon_1+\varepsilon_2[$$

etc (see the process tree in Diagram B).
    \end{itemize}

More generally, let us consider for all $n\geq0$ the index set $I_n=\{0,1,...,n\} $ and the index sets
\begin{equation}
I_n^+=\{i\in I_n\ / \sigma_i=+\}\qquad \hbox{and}\qquad I^-_n=\{i\in I_n\ / \sigma_i=-\}.
  \end{equation}
We have $I_n=I_n^+\cup I^-_n$ for all $n\geq 0$.

For all $n\geq0$ and for all $j_n=\sigma_0...\sigma_n\in \Lambda_n$, where $\Lambda_n$ is the index set defined by (\ref{DefSet}), the stretching process gives the following intervals at the $step(n)$:
\begin{equation}
X_n^{j_n}=\Big]a-\sum_{i\in I_n^-}\varepsilon_i\ ,\ b+\sum_{i\in I_n^+}\varepsilon_i\Big[,
\end{equation}
and by stretching left and right the open interval $X_n^{j_n}$ we obtain at the $step(n+1)$ the open intervals
\begin{equation}
X_{n+1}^{j_n-}=\Big]a-\sum_{i\in I_n^-}\varepsilon_i-\varepsilon_{n+1}\ ,\ b+\sum_{i\in I_n^+}\varepsilon_i\Big[,\quad\hbox{and}\quad X_{n+1}^{j_n+}=\Big]a-\sum_{i\in I_n^-}\varepsilon_i\ ,\ b+\sum_{i\in I_n^+}\varepsilon_i+\varepsilon_{n+1}\Big[.
\end{equation}
\gesp

\unitlength=0.85cm
\begin{picture}(11,7)
\put(7.7,7.3){$]\ a\ ,\ b\ [$}
\put(6.8,5.7){$ X_0^-$}

\put(9.4,5.7){$X_0^+$}


\put(4.8,3.5){$X_1^{--}$}

\put(6.8,3.5){$X_1^{-+}$}

\put(9.2,3.5){$X_1^{+-}$}

\put(11.2,3.5){$X_1^{++}$}


\put(3.95,1){$X_2^{---}$}

\put(5.5,-0.3){$X_2^{--+}$}

\put(6,1){$X_2^{-+-}$}

\put(7.5,-0.3){$X_2^{-++}$}

\put(8.3,1){$X_2^{+--}$}

\put(9.8,-0.3){$X_2^{+-+}$}

\put(10.45,1){$X_2^{++-}$}

\put(12,-0.3){$X_2^{+++}$}

\put(8.4,7.1){\vector(-1,-1){1}}

\put(8.4,7.1){\vector(1,-1){1}}


\put(6.8,5.5){\vector(-1,-1){1.6}}

\put(7.1,5.5){\vector(0,-3){1.5}}

\put(9.4,5.5){\vector(0,-3){1.5}}

\put(9.8,5.5){\vector(1,-1){1.6}}


\put(5.1,3.2){\vector(-1,-3){.6}}

\put(5.1,3.2){\vector(1,-3){1}}


\put(7.1,3.2){\vector(-1,-3){.6}}

\put(7.1,3.2){\vector(1,-3){1}}


\put(9.4,3.2){\vector(-1,-3){.6}}

\put(9.4,3.2){\vector(1,-3){1}}


\put(11.5,3.2){\vector(-1,-3){.6}}

\put(11.5,3.2){\vector(1,-3){1}}

\put(0.1,-1.5){\footnotesize {\it Diagram B: Tree that illustrates the construction of an expanding fractal family of topological spaces}}
\put(1.4,-2){\footnotesize {\it\qquad  via left stretching and right stretching process of the open interval $]a,b[$.}}

\thicklines
\end{picture}
\vskip 2.5cm

For all $n\geq0$, and for all $j_n\in \Lambda_n$, the open interval $X_n^{j_n}$ of $\rR$ is endowed with the topology ${\cal F}_n^{j_n}$ induced by the usual topology of $\rR$. It is not difficult to prove that this family of topological spaces $\Big(X_n^{j_n}, {\cal F}_n^{j_n}\Big)_{{{j_n \in\Lambda_n}\atop n\geq0}}$ verifies the following:

\begin{enumerate}
  \item $Card\ \Lambda_n<Card\ \Lambda_{n+1}$ for all $n\geq0$;
  \item $(X_n^{j_n}, {\cal F}_n^{j_n})$ is a topological space for all $n\geq0$ and for all $j_n\in \Lambda_n$;
  \item For a given $n\geq0$, the topologies ${\cal F}_n^{j_n}$ are equivalent for all $j_n\in \Lambda_n$ since the $X_n^{j_n}$, $j_n\in I_n$ are open intervals of $\rR$;
  \item For all $n\geq0$ and for all $j_{n+1}=\sigma_0...\sigma_n\sigma_{n+1}\in \Lambda_{n+1}$, there exists a unique $j_n\in \Lambda_n$ given by $j_n=\sigma_0...\sigma_n$ such that $X_n^{j_n}\subset X_{n+1}^{j_{n+1}}$ and ${\cal F}_n^{j_n}=\{O\cap X_n^{j_n}\ /\ O\in X_{n+1}^{j_{n+1}}\}$;
  \item For all $n\geq0$ and for all $j_{n}=\sigma_0...\sigma_n\in \Lambda_{n}$, there exist two indexes $j_{n+1}\in I_{n+1}$ given by $j_{n+1}=\sigma_0...\sigma_n+$ and $j_{n+1}=\sigma_0...\sigma_n-$ such that ${\cal F}_n^{j_n}\subset{\cal F}_{n+1}^{j_{n+1}}$ and ${\cal F}_n^{j_n}=\{O\cap X_n^{j_n}\ /\ O\in X_{n+1}^{j_{n+1}}\}$.
\end{enumerate}
Consequently the family $\Big(X_n^{j_n}, {\cal F}_n^{j_n}\Big)_{{{j_n \in\Lambda_n}\atop n\geq0}}$ is a fractal family of topological spaces. If we define the coproduct $X_n=\coprod_{j_n\in \Lambda_n}X_n^{j_n}=\bigcup_{j_n\in \Lambda_n}X_n^{j_n}\times \{j_n\}$ together with the coproduct topology ${\cal F}_n$, then by Theorem \ref{P5}, the family $\Big(X_n, {\cal F}_n\Big)_{n\geq0}$ is an expanding family of topological spaces, and by Theorem \ref{P1}, ${\cal F}_n\subset{\cal F}_{n+1}$. Moreover the topological space $]a,b[$ is a locally expandable topological space, indeed
for all $n\geq0$ and for all $j_n\in \Lambda_n$ the open interval $]a,b[$ is homeomorphic to the open interval $X_n^{j_n}$, then $]a,b[$ is locally homeomorphic to $X_n^{j_n}$ for all $j_n\in \Lambda_n$ and then there exists a local multi-homeomorphism between $]a,b[$ and $X_n=\di\coprod_{j_n\in \Lambda_n}X_n^{j_n}$, which yields by Definition \ref{D6} that $]a,b[$ is locally expandable. Therefore by Theorem \ref{Th3} $]a,b[$ has a topological expansion. We have the following proposition:

\begin{prop}\label{P7}
Any finite open interval of $\rR$ is locally expandable by the stretching process on $\rR$.
\end{prop}

\ni{\bf Proof.} The proof is straight forward from the description of the previous stretching process.

\rightline\Box

\begin{figure}[h]
\begin{center}
\includegraphics[width=10cm]{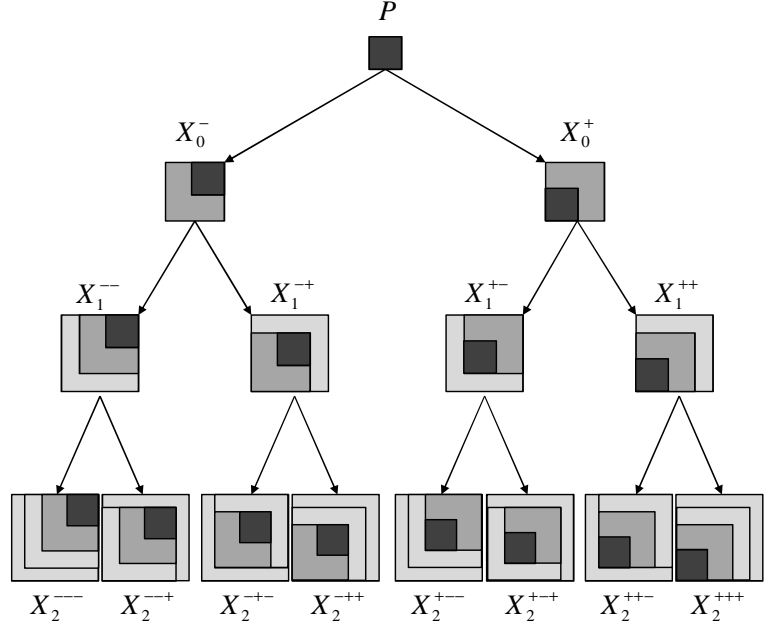}
\centering\caption{ {\footnotesize This tree illustrates the fractal constituent spaces of the coproduct $(X_n, {\cal F}_n)$ for $n=0,1,2$. The existence of a fractal family of topological spaces  $\{(X_n^{j_n}, {\cal F}_n^{j_n}),/ j_n \in\Lambda_n\ \hbox{and}\  0\leq n\leq2\}$ with $Card \Lambda_n=2^{n+1}$, makes any open box in $\rR^2$ locally expandable. A stretching process on $\rR^2$ is obtained by left-down stretching and right-up stretching of the box $P$ in the plane, which corresponds to the left and right stretching of each side of the open box $P$ from step(0) to step(2).
}}\label{Fig1}
\end{center}
\end{figure}

\begin{rem}
The topological expansion of $]a,b[$ is not unique since there exists an infinite number of sequences $(\varepsilon_n)_{n\geq0}$ such that $\varepsilon_n<...<\varepsilon_1<\varepsilon_0$ for all $n\geq 0$, which gives an infinite number of stretching processes.
\end{rem}

\begin{rem}
The open interval $\ ]a,b[\subset\rR\ $ is not expandable by the stretching process if $a=-\infty$ or $b=+\infty$.
In particular $\rR$ is not expandable by the stretching process.
\end{rem}

\begin{description}
  \item[Stretching process on $\rR^n$:] the stretching process on $\rR$ can be generalized to $\rR^n$.
\end{description}

Let us consider $$P=\Big\{(x_1,...,x_2)\in \rR^n\ /\ a_k<x_k<b_k,\ \ \forall k=1,...,n\Big\}$$

\begin{figure}[h!]
\begin{center}
\includegraphics[width=10cm]{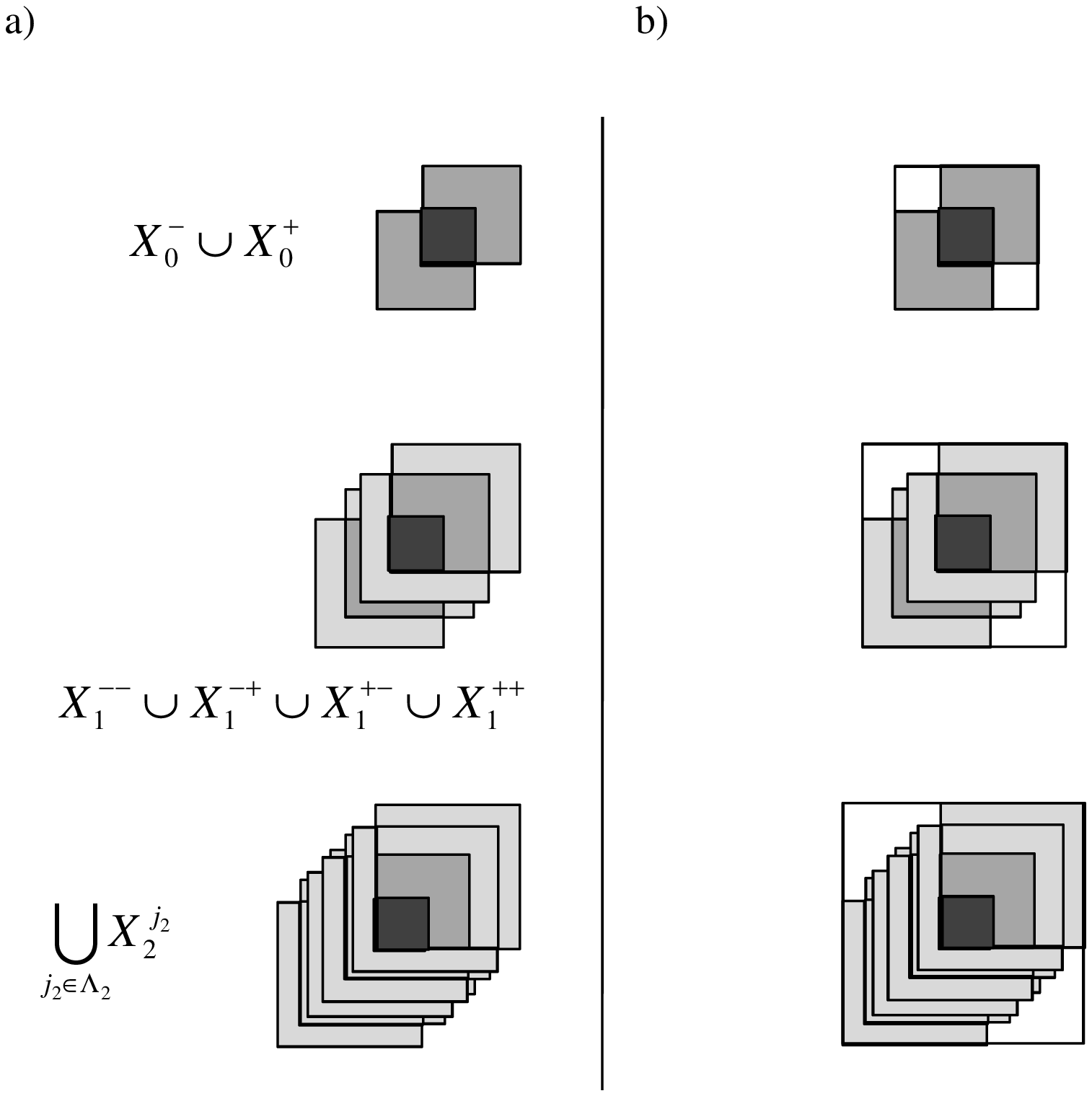}
\centering\caption{ {\footnotesize The non disjoint union of all constituents of the stretching process illustrated in Fig.1 from $step(0)$ to $step(2)$ gives different areas that index the expansion states of the square box in a). The union $X_2^{---}\cup X_2^{--+}\cup X_2^{-+-}\cup X_2^{-++}\cup X_2^{+--}\cup X_2^{+-+}\cup X_2^{++-}\cup X_2^{+++}=\cup_{j_2\in\Lambda_2}X_n^{j_2}$ for example is the new surface of the square box obtained at the $step(2)$ after left and right stretching process. This new surface is compared to the surface of bigger square box that delimitate it in b), where the white surface within the big square box is missing in the stretching process of the box $P$, which makes the fractal dimension of the obtained stretched surface less that 2 since it is not a surface of full square. .
}}\label{Fig2}
\end{center}
\end{figure}
\ni an open box of $\rR^n$, where $a_k$ and $b_k$ are finite real numbers for $k=1,...,n$, and let $(\varepsilon_n)_{n\geq0}$ be a strictly decreasing sequence of real numbers. By applying the stretching process on $\rR$ to each open interval $]a_k,b_k[$ for k=1,...,n, we obtain a fractal family of topological spaces $\Big(X_n^{j_n}, {\cal F}_n^{j_n}\Big)_{{{j_n \in\Lambda_n}\atop n\geq0}}$ given for all $n\geq0$ and for all $j_n\in \Lambda_n$ by
\begin{equation}
X_n^{j_n}=\Big\{(x_1,...,x_2)\in \rR^n\ /\ x_k\in \Big]a_k-\sum_{i\in I_n^-}\varepsilon_i\ ,\ b_k+\sum_{i\in I_n^+}\varepsilon_i\Big[,\quad \forall k=1,...,n\ \Big\}
\end{equation}
endowed with the topology ${\cal F}_n^{j_n}$ induced by the usual topology of $\rR^n$.
The coproduct associated to the fractal family $\Big(X_n^{j_n}, {\cal F}_n^{j_n}\Big)_{{{j_n \in\Lambda_n}\atop n\geq0}}$ for all $n\geq0$ and for all $j_n\in \Lambda_n$ is given by
\begin{equation}
X_n=\bigcup_{j_n\in \Lambda_n}\Big(\prod_{k=1}^n\Big]a_k-\sum_{i\in I_n^-}\varepsilon_i\ ,\ b_k+\sum_{i\in I_n^+}\varepsilon_i\Big[\Big)\times\{j_n\}
\end{equation}
By Proposition \ref{P7}, the open box $P$ is locally expandable (see illustration of a locally expandable box in $\rR^2$, Fig.1) and we have the following proposition:

\begin{prop}
Any open box of $\rR^n$ is locally expandable by the stretching process on $\rR^n$.
\end{prop}

\ni{\bf Proof.}
The proof is straight forward from the construction.

\rightline\Box

\begin{rem}
An open box of $\rR^2$ is expandable, where the stretching process from $step(0)$ to $step(2)$ is illustrated in Fig.1. The non disjoint union of the obtained fractal family $\{(X_n^{j_n}, {\cal F}_n^{j_n}),/ j_n \in\Lambda_n\ \hbox{and}\  0\leq n\leq2\}$ is illustrated in Fig.2, and its is not difficult to see that the fractal dimension of the surface $\bigcup_{j_n\in\Lambda_n}X_n^{j_n}$ for $0\leq n\leq2$ is strictly less than 2.
\end{rem}

\subsection{Fundamental Application: the Fractal Manifold}

\subsubsection{Coproduct Topology associated to the Fractal Manifold}

Let ${\cal M}$ be a fractal manifold. By Theorem \ref{Th2}, ${\cal M}$ is locally homeomorphic via a family of local homeomorphisms
to the fractal family of topological spaces $\Big(\bigcup _{{\delta_0}...{\delta_{n}}} N_{\delta_0...\delta_{n}}^{j_n},\ {\cal T}_{n}^{j_n}\Big)_{{j_n\in \Lambda_n\atop n\geq0}}$. Let us consider the coproduct family  $\Big(N_n,{\cal T}_n\Big)_{n\geq0}$ associated to $\Big(\bigcup _{{\delta_0}...{\delta_{n}}} N_{\delta_0...\delta_{n}}^{j_n},\ {\cal T}_{n}^{j_n}\Big)_{{j_n\in \Lambda_n\atop n\geq0}}$ given for all $n\geq0$ by
\begin{equation}\label{FO2}
\begin{array}{lll}
N_n & =\di \coprod_{j_n\in\Lambda_n}\Big(\bigcup _{{\delta_0}...{\delta_{n}}} N_{\delta_0...\delta_{n}}^{j_n}\Big)=\bigcup_{j_n\in\Lambda_n}\Big(\bigcup _{{\delta_0}...{\delta_{n}}} N_{\delta_0...\delta_{n}}^{j_n}\Big)\times\{j_n\}\\
{\cal T}_n & =\di \Big\{ \bigcup_{j_n\in\Lambda_n}\Omega_n^{j_n}\times\{j_n\}\ /\ \forall j_n\in\Lambda_n,\ \Omega_n^{j_n}\in{\cal T}_n^{j_n}\ \Big\}.
\end{array}
\end{equation}

\begin{prop}\label{P6}
The coproduct family $\Big(N_n, {\cal T}_n\Big)_{n\geq0}$ is an expanding family of topological spaces.
\end{prop}

\ni{\bf Proof.} By Theorem \ref{Th5}, the family $\Big(\bigcup _{{\delta_0}...{\delta_{n}}} N_{\delta_0...\delta_{n}}^{j_n},\ {\cal T}_{n}^{j_n}\Big)_
{{j_n\in\Lambda_n\atop n\geq0}}$ is a fractal family of topological spaces, then  by Theorem \ref{P5} $\Big(N_n, {\cal T}_n\Big)_{n\geq0}$ is an expanding family of topological spaces.

\rightline\Box

\begin{prop}
For all ${n\geq0}$, $\ {\cal T}_n\subset {\cal T}_{n+1}$.
\end{prop}

\ni{\bf Proof.} The proof is a direct consequence of Theorem \ref{P1}.

\rightline\Box

\subsubsection{Topological Expansion of the Fractal Manifold}

From the Theorem \ref{Th1}, Diagram A, one can understand the expansion of the fractal manifold $\cal M$ in the local coordinates system at each step. Indeed, for all $P\in\cal M$ there exists a local multi-chart at $P$ given by
\begin{equation}
\Big(\Omega,\psi_1,\psi'_1\Big)\qquad \hbox{for the}\quad step(0)
\end{equation}
\begin{equation}
\Big(\Omega,\psi_2,\psi'_2,\psi_3,\psi'_3\Big)\qquad \hbox{for the}\quad step(1)
\end{equation}
\begin{equation}
\Big(\Omega,\psi_4,\psi'_4,\psi_5,\psi'_5,\psi_6,\psi'_6,\psi_7,\psi'_7\Big)\qquad \hbox{for the}\quad step(2),
\end{equation}
etc,
\begin{equation}\label{FO8}
\Big(\Omega,\psi_{2^n},\psi'_{2^n},\psi_{2^n+1},\psi'_{2^n+1},\ldots\ldots,\psi_{2^{n+1}-1},\psi'_{2^{n+1}-1}\Big)\qquad \hbox{for the}\quad step(n),
\end{equation}
where $\psi_1=\varphi_1$ and $\psi'_1=T_1\circ\varphi_1$, $\psi_2=\varphi_2\circ\varphi_1$, and $\psi'_2=T_2\circ\varphi_2\circ\varphi_1$, $\psi_3=\varphi_3\circ T_1\circ\varphi_1$ and $\psi'_3=T_3\circ\varphi_3\circ T_1\circ\varphi_1$, etc.

If we use the index set $\Lambda_n$ given by (\ref{DefSet}) with $Card\ \Lambda_n=2^{n+1}$, then we will have a family $(\psi_n^{j_n})_{j_n\in\Lambda_n}$ that represents the $2^{n+1}$ local homeomorphisms used in the local multi-chart (\ref{FO8}) at the $step(n)$ for all $n\geq0$, therefore the local multi-chart (\ref{FO8}) at the $step(n)$ can be denoted as follow $\Big(\Omega,(\psi_n^{j_n})_{j_n\in\Lambda_n}\Big)$, where for $n=0$ this chart is a triplet, for $n=1$ it is a quintuplet, etc.
Since for all $n\geq0$ and for all $j_n\in \Lambda_n$ we have
\begin{equation}
\psi_n^{j_n}:\ \Omega\longrightarrow\psi_n^{j_n}(\Omega)\subset\bigcup_{\delta_0...\delta_n}N_{\delta_0...\delta_n}^{j_n},
\end{equation}
then for all $n\geq0$, the disjoint union $\di\coprod_{j_n\in\Lambda_n}\psi_n^{j_n}(\Omega)=\bigcup_{j_n\in\Lambda_n}\psi_n^{j_n}(\Omega)\times\{j_n\}$ belongs to the topology ${\cal T}_n$ given by (\ref{FO2}).
The open set $\di\coprod_{j_n\in\Lambda_n}\psi_n^{j_n}(\Omega)$ is a neighborhood of $\bigcup_{j_n\in\Lambda_n}\psi_n^{j_n}(P)$ in the topological space $(N_n, {\cal T}_n)$, where $\di\bigcup_{j_n\in\Lambda_n}\psi_n^{j_n}(P)$ is the local representation of $P\in {\cal M}$ at the $step(n)$ in the local multi-chart $\Big(\Omega, (\psi_n^{j_n})_{j_n\in\Lambda_n}\Big)$. The number of local coordinates of $P$ in the local multi-chart $\Big(\Omega, (\psi_n^{j_n})_{j_n\in\Lambda_n}\Big)$ is increasing from one step to another. We can express the expansion of the fractal manifold $\cal M$  via its local expansion everywhere:

\begin{prop}\label{P8}
If $\cal M$ is a fractal manifold, then $\cal M$ is locally multi-homeomorphic to $(N_n, {\cal T}_n)$ for all $n\geq0$.
\end{prop}

\ni{\bf Proof.}
Let $\cal M$ be a fractal manifold. For all $P\in {\cal M}$ and $n\geq0$, there exists a local chart $\Big(\Omega, (\psi_n^{j_n})_{j_n\in\Lambda_n}\Big)$ such that
$\Omega$ is a neighborhood of $P$ in ${\cal M}$, $\di\coprod_{j_n\in \Lambda_n}\psi_n^{j_n}(\Omega)\in {\cal T}_n$ and
for all $j_n\in \Lambda_n$, $\psi_n^{j_n}: \Omega\longrightarrow \psi_n^{j_n}(\Omega)$ is an homeomorphism.
Therefore by Definition \ref{D3}, $\cal M$ is locally multi-homeomorphic to $(N_n, {\cal T}_n)$ for all $n\geq0$, which completes the proof.

\rightline\Box

\begin{cor}\label{Cor1}
If $\cal M$ is a fractal manifold, then $\cal M$ is locally expandable.
\end{cor}

\ni{\bf Proof.}
Let $\cal M$ be a fractal manifold, then there exists a fractal family of topological spaces given by $\Big(\di\bigcup _{{\delta_0}...{\delta_{n}}} N_{\delta_0...\delta_{n}}^{j_n},\ {\cal T}_{n}^{j_n}\Big)_{{j_n\in \Lambda_n\atop n\geq0}}$ such that  $\cal M$ is locally multi-homeomorphic to the associated coproduct topological space $(N_n, {\cal T}_n)$ for all $n\geq0$ as proved in Proposition \ref{P8}, then by Definition \ref{D6} we obtain that ${\cal M}$ is locally expandable, which completes the proof.

\rightline\Box

\begin{thm}\label{Th7}
If $\cal M$ is a fractal manifold, then $\cal M$ has a topological expansion.
\end{thm}

\ni{\bf Proof.} Since $\cal M$ is a fractal manifold, then by Corollary \ref{Cor1}, ${\cal M}$ is locally expandable, and then by Theorem \ref{Th3}, ${\cal M}$ has a topological expansion.

\rightline\Box

\begin{rem}
The local transformations of the fractal manifold ${\cal M}$  from one step to another are locally identified via local multi-homeomorphisms to the local transformations of the expanding topological space $(N_0, {\cal T}_0)$. Therefore the more ${\cal M}$ expands, the more ${\cal M}$ is locally multi-homeomorphic to a finer topological space.
\end{rem}

\section{Overview and Conclusion}

In this framework we have established a new formalism that leads:

i) to define properly a topological expansion of a topological space via disjoint union topology of family of image of topological spaces $X_n^{j_n}$ for all $n\geq0$, where the image of $X_n^{j_n}$ is given by the product $X_n^{j_n}\times \{j_n\}$;

ii) to make possible the comparison of two different topological spaces $(X,\tau_1)$ and $(Y,\tau_2)$ when $X$ is a subspace of $Y$;

iii) to quantify a topological expansion in different indexed expansion states, and thanks to ii), to end up with a clear causality between space expansion and the finer topology of the indexed expansion states: the more a topological space expands, the finer the topology of its indexed states is.

The main lead in this new formalism is the fractal manifold model. Indeed, contrary to the Riemannian or pseudo-Riemannian manifold used in general relativity, the fractal manifold model presents two fundamental aspects: a variable geometry \cite{BF3}, and a variable topology \cite{HP}. The coexistence of these two aspects leads to understand the space-time dynamic via topological variation of the space as it was elaborated within this work.

It is known that the topology is invariant under continuous deformations of the space, including space expansion, and then it is quite impossible to find criteria using classical topology that give a new understanding of the space expansion via topology (or to understand what happens to the topology when the space expands). However, this new formalism makes the fractal topology sensitive to the continuous deformations of the space such as its expansion, and allows to find out a causality between space expansion and topology variation. The reason behind that is the main structure of the fractal manifold \cite{BF1}. Indeed the continuous deformation of the space on fractal manifold is constructed via a discontinuous process: a double family of local homeomorphisms from one step to another that double the local properties of the topological space and generates local expansion of the space. This procedure in describing the continuous deformation of the space such as its expansion is actually made by quantifying splits of the local properties of the space. This sequence of splits makes the fractal topology on fractal manifold locally sensitive to the continuous deformations of the space, which gives the causality between space expansion and variation of the fractal topology.

The global expansion of the space, within this framework, can be understood via its local expansion everywhere that preserves the homogeneity of the space, which is one of the fundamental principle in cosmology. This framework provides a new insight and define causalities related to the variation of topology as the space-time expands, which could open a new window on topology applications in cosmology.


\bibliographystyle{spmpsci}      

\end{document}